\documentclass[leqno]{jktr} \def\trimmarks{}
\def\draftnote{\today\quad\currenttime\qquad
  revised version \qquad\jobname}

\usepackage{amsmath,mathrsfs,graphicx,subfigure}
\usepackage{mathptmx,courier}
\usepackage[scaled=0.95]{helvet}

\title{\large Equivalence of symmetric union diagrams}

\author{Michael Eisermann}
\address{Institut Fourier, Universit\'e Grenoble I, France \\[-2pt]
Michael.Eisermann@ujf-grenoble.fr}
\author{Christoph Lamm}
\address{Karlstra{\ss}e 33, 65185 Wiesbaden, Germany \\[-2pt]
Christoph.Lamm@web.de}
\newcommand{\qed}{\hfill$\Box$}

\newcommand{\Z}{\mathbb{Z}}
\newcommand{\R}{\mathbb{R}}

\newcommand{\into}{\hookrightarrow}

\newcommand{\isoto}{\mathrel{\xrightarrow{_\sim}}} 
\newcommand{\cs}{\mathbin{\sharp}}

\newsavebox{\boxaocr}\savebox{\boxaocr}%
{\smash{\raisebox{-0.65ex}{\includegraphics[height=2.5ex]{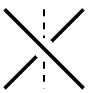}}}}
\newsavebox{\boxaucr}\savebox{\boxaucr}%
{\smash{\raisebox{-0.65ex}{\includegraphics[height=2.5ex]{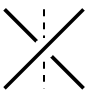}}}}
\newsavebox{\boxahcr}\savebox{\boxahcr}%
{\smash{\raisebox{-0.65ex}{\includegraphics[height=2.5ex]{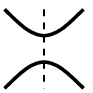}}}}
\newsavebox{\boxavcr}\savebox{\boxavcr}%
{\smash{\raisebox{-0.65ex}{\includegraphics[height=2.5ex]{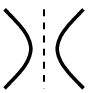}}}}

\newcommand{\aocr}{\usebox{\boxaocr} }
\newcommand{\aucr}{\usebox{\boxaucr} }
\newcommand{\ahcr}{\usebox{\boxahcr} }
\newcommand{\avcr}{\usebox{\boxavcr} }

\begin{document} %%%%%%%%%%%%%%%%%%%%%%%%%%%%%%%%%%%%%%%%%%%%%%%%%%%%%%%%%%%%
%%%%%%%%%%%%%%%%%%%%%%%%%%%%%%%%%%%%%%%%%%%%%%%%%%%%%%%%%%%%%%%%%%%%%%%%%%%%%

\maketitle

% \headline{90mm}{\eprintinfo}\vspace*{-7mm}

\begin{abstract}
  Motivated by the study of ribbon knots we explore symmetric unions, 
  a beautiful construction introduced by Kinoshita and Terasaka 
  50 years ago.  It is easy to see that every symmetric union 
  represents a ribbon knot, but the converse is still an open problem.  
  Besides existence it is natural to consider the question of uniqueness.
  In order to attack this question we extend the usual Reidemeister moves 
  to a family of moves respecting the symmetry, and consider 
  the symmetric equivalence thus generated.  This notion being 
  in place, we discuss several situations in which a knot 
  can have essentially distinct symmetric union representations.
  We exhibit an infinite family of ribbon two-bridge knots 
  each of which allows two different symmetric union representations.
\end{abstract}
  


\medskip

\begin{center} \it
  Dedicated to Louis H.\ Kauffman 
  on the occasion of his 60th birthday
\end{center}

%%%%%%%%%%%%%%%%%%%%%%%%%%%%%%%%%%%%%%%%%%%%%%%%%%%%%%%%%%%%%%%%%%%%%%%%%%%%%

\section{Motivation and background}

Given a ribbon knot $K$, Louis Kauffman emphasized in his 
course notes \textit{On knots} \cite[p.\ 214]{Kauffman:OnKnots} 
that ``in some algebraic sense $K$ looks like a connected 
sum with a mirror image.  Investigate this concept.''
Symmetric unions are a promising geometric counterpart of this analogy, 
and in continuation of Kauffman's advice, their investigation shall be  
advertised here.  Algebraic properties, based on a refinement of the 
bracket polynomial, will be the subject of a forthcoming paper.
\hfill 
\textit{Happy Birthday, Lou!}

\subsection{Symmetric unions}

\begin{figure}[hbtp]
  \centering
  \includegraphics[scale=0.6]{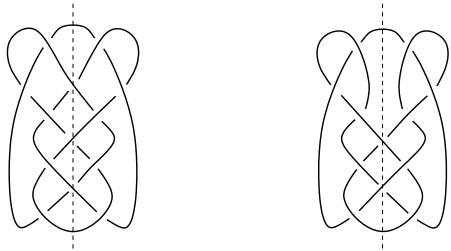}
  \caption{Two symmetric union presentations of the knot $9_{27}$}
  \label{fig:Knot-9-27}
\end{figure}

In this article we consider symmetric knot diagrams and study 
the equivalence relation generated by symmetric Reidemeister moves.
Figure \ref{fig:Knot-9-27} shows two typical examples of such diagrams.
Notice that we allow any number of crossings on the axis --- they necessarily 
break the mirror symmetry, but this defect only concerns the crossing sign
and is localized on the axis alone.

We are particularly interested in \emph{symmetric unions}, where %% for a \emph{union} 
we require the diagram to represent a knot (that is, a one-component link) 
that traverses the axis in exactly two points that are not crossings.
In other words, a symmetric union looks like the connected sum of 
a knot $K_+$ and its mirror image $K_-$, with additional crossings 
inserted on the symmetry axis.  

Conversely, given a symmetric union, one can easily recover the two partial knots
$K_-$ and $K_+$ as follows:  they are the knots on the left and on the right 
of the axis, respectively, obtained by cutting open each crossing on the axis,
according to $\aocr \mapsto \avcr$ or $\aucr \mapsto \avcr$.
The result is a connected sum, which can then be split by one final cut 
$\ahcr \mapsto \avcr$ to obtain the knots $K_+$ and $K_-$, as desired.
(In Figure \ref{fig:Knot-9-27}, for example, we find the partial knot $5_2$.)
%% and its mirror image.)

\subsection{Ribbon knots}

\begin{figure}[hbtp]
  \centering
  \includegraphics[scale=0.8]{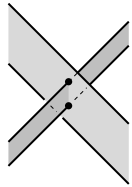}
  \caption{An immersed disk with ribbon singularity}
  \label{fig:RibbonSingularity}
\end{figure}

\begin{figure}[hbtp]
  \centering
  \subfigure[$8_{20}$]{\includegraphics[scale=0.8]{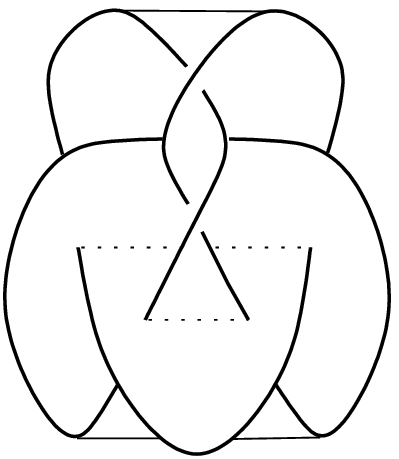}} \hspace{20mm}
  \subfigure[$10_{87}$]{\includegraphics[scale=0.8]{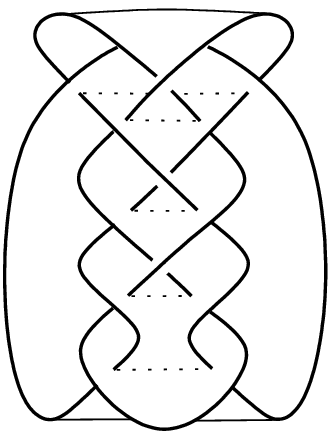}}
  \caption{The knots $8_{20}$ and $10_{87}$ represented as symmetric unions. 
    The figure indicates the resulting symmetric ribbon with twists.}
  \label{fig:SymmetricRibbon}
\end{figure}

Symmetric unions have been introduced in 1957 by 
Kinoshita and Terasaka \cite{KinoshitaTerasaka:1957}.  
Apart from their striking aesthetic appeal, symmetric 
unions appear naturally in the study of ribbon knots.
We recall that a knot $K \subset \R^3$ is a \emph{ribbon knot} 
if it bounds a smoothly immersed disk $\mathbb{D}^2 \looparrowright \R^3$
whose only singularities are ribbon singularities 
as shown in Figure \ref{fig:RibbonSingularity}: 
two sheets intersecting in an arc whose preimage consists of a properly
embedded arc in $\mathbb{D}^2$ and an embedded arc interior to $\mathbb{D}^2$.
Figure \ref{fig:SymmetricRibbon} displays two examples.

Put another way, a knot $K \subset \R^3$ is a ribbon knot if and only if it bounds 
a locally flat disk $\mathbb{D}^2 \into \R^4_+ = \{ x \in  \R^4 \mid x_4 \ge 0 \}$
without local minima.  More generally, if $K$ bounds an arbitrary 
locally flat disk in $\R^4_+$, then $K$ is called a \emph{slice knot}.
It is a difficult open question whether every smoothly slice knot 
is a ribbon knot.  For a general reference see \cite{Livingston:2005}.

For the rest of this article we will exclusively work in the smooth category.

\subsection{Which ribbon knots are symmetric unions?}

While it is easy to see that every symmetric union represents a ribbon knot, 
as in Figure \ref{fig:SymmetricRibbon}, the converse question is still open.  
The following partial answers are known:
\begin{itemize}
\item
  There are $21$ non-trivial prime ribbon knots with at most $10$ crossings.
  By patiently compiling an exhaustive list, Lamm \cite{Lamm:2000,Lamm:Talk} 
  has shown that each of them can be presented as a symmetric union.
  %% For the convenience of the reader, and as a decorative illustration,
  %% we reproduce these diagrams in Figure \ref{fig:Compilation}.
\item
  In 1975, Casson and Gordon \cite{CassonGordon:1975} 
  exhibited three infinite families of two-bridge ribbon knots,
  and Lamm \cite{Lamm:Talk} has shown that each 
  of them can be presented as a symmetric union.
  Recently, Lisca \cite{Lisca:2007} has shown that the
  three Casson-Gordon families exhaust all two-bridge ribbon knots.
\end{itemize}

\begin{remark}
  Presenting a given knot $K$ as a symmetric union
  is one way of proving that $K$ is a ribbon knot,
  and usually a rather efficient one, too.  
  The explicit constructions presented here 
  have mainly been a matter of patience, and it is fair 
  to say that symmetry is a good guiding principle.
\end{remark}

\begin{figure}[hbtp]
  \centering
  \includegraphics[scale=0.55]{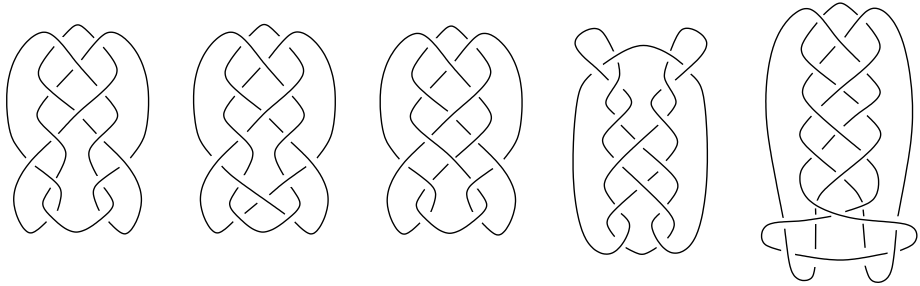}
  \caption{Symmetric unions representing %% the knots 
    11a28, 11a35, 11a36, 11a96, and 11a164.  This proves 
    that they are ribbon, hence smoothly slice.}
  \label{fig:Knots-11x}
\end{figure}

\begin{figure}[hbtp]
  \centering
  \includegraphics[scale=0.65]{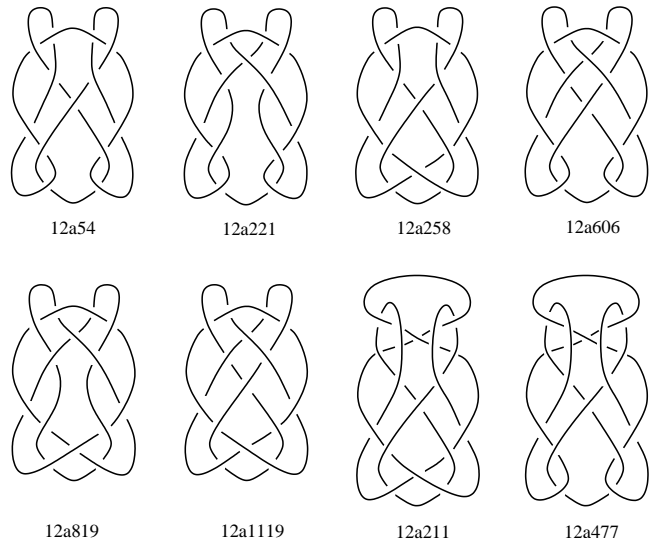}
  \caption{Symmetric union presentations for some ribbon knots 
    with 12 crossings, all with partial knot $C(2,1,1,2)$.}
  \label{fig:Knots-12x}
\end{figure}

\begin{example}  
  When venturing to knots with $11$ crossings, there still remain, 
  at the time of writing, several knots that are possibly ribbon
  in the sense that their algebraic invariants do not obstruct this.
  It remains to explicitly construct a ribbon --- 
  or to refute this possibility by some refined argument.
  According to Cha and Livingston \cite{ChaLivingston:2003v6}, as of March 2006, 
  there remained eleven knots of which it was not known whether 
  they were slice.  Figure \ref{fig:Knots-11x} solves this question 
  for five of them by presenting them as symmetric unions.  
  In the same vein, Figure \ref{fig:Knots-12x} displays some 
  $12$-crossing knots (which all have identical partial knots).
\end{example}

\pagebreak[5]

\begin{remark}[notation]
  For knots with up to $10$ crossings we follow the traditional numbering 
  of Rolfsen's tables \cite{Rolfsen:1990} with the correction by Perko \cite{Perko:1974}.
  For knots with crossing number between $11$ and $16$ we use 
  the numbering of the KnotScape library \cite{Knotscape}.
  Finally, $C(a_1,a_2,\dots,a_n)$ is Conway's notation 
  for two-bridge knots, see \cite[\textsection2.1]{Kawauchi:1996}.
\end{remark}

\begin{question}
  Can every ribbon knot be presented as a symmetric union?
  This would be very nice, but practical experience suggests 
  that it is rather unlikely.  A general construction, 
  if it exists, must be very intricate.
\end{question}

\begin{remark}
  The search for symmetric union presentations can be automated,
  and would constitute an interesting project at the intersection
  of computer science and knot theory.  The idea is to produce
  symmetric union diagrams in a systematic yet efficient way,
  and then to apply KnotScape to identify the resulting knot.
  Roughly speaking, the first step is easy but the second usually 
  takes a short while and could turn out to be too time-consuming.  
  Library look-up should thus be used with care, and the production 
  of candidates should avoid duplications as efficiently as possible.
  (This is the non-trivial part of the programming project.)

  In this way one could hope to find symmetric union 
  presentations for all remaining $11$-crossing knots, 
  and for many knots with higher crossing numbers as well.  
  Such a census will yield further evidence of how large 
  the family of symmetric unions is within the class of ribbon knots 
  --- and possibly exhibit ribbon knots that defy symmetrization.
  No such examples are known at the time of writing.  
  Of course, once a candidate is at hand, a suitable obstruction
  has to be identified in order to prove that it cannot be represented 
  as a symmetric union.  (This is the non-trivial mathematical part.)
\end{remark}

\subsection{Symmetric equivalence}

Besides the problem of \emph{existence} it is natural to consider 
the question of \emph{uniqueness} of symmetric union representations.
Motivated by the task of tabulating symmetric union diagrams 
for ribbon knots, we are led to ask when two such diagrams 
should be regarded as equivalent.  
One way to answer this question is to extend the usual 
Reidemeister moves to a family of moves respecting the symmetry, 
as explained in \textsection\ref{sub:SymmetricMoves}.

\begin{figure}[hbtp]
  \centering
  \includegraphics[scale=0.55]{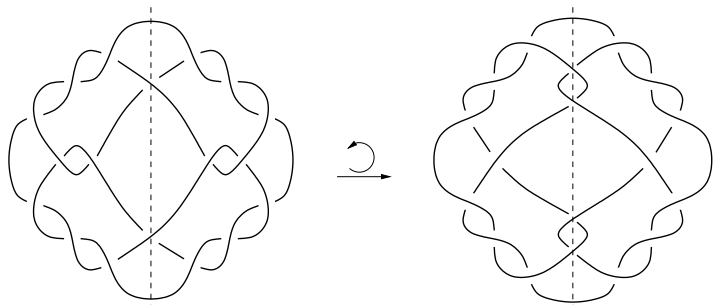}
  \caption{A symmetric union with two symmetry axes}
  \label{fig:TwoAxes}
\end{figure}

\begin{example}
  It may well be that two symmetric union representations are equivalent 
  (via the usual Reidemeister moves), but that such an equivalence is not symmetric,
  that is, the transformation cannot be performed in a symmetric way.
  One possible cause for this phenomenon is the existence of 
  \emph{two} axes of symmetry.  The simplest (prime) example of this type 
  seems to be the knot 16n524794 shown in Figure \ref{fig:TwoAxes}:
  the symmetric unions have partial knots $6_1$ and $8_{20}$, 
  respectively, and thus cannot be symmetrically equivalent.
\end{example}

\begin{remark}
  For a symmetric union representing a knot $K$ 
  with partial knots $K_+$ and $K_-$ the determinant
  satisfies the product formula $\det(K) = \det(K_+) \det(K_-)$
  and is thus a square.  This was already noticed by 
  Kinoshita and Terasaka \cite{KinoshitaTerasaka:1957} 
  in the special case that they considered;  
  for the general case see \cite{Lamm:2000}.
  For a symmetric union with two symmetry axes this means 
  that the determinant is necessarily a fourth power.
\end{remark}

\begin{figure}[hbtp]
  \centering
  \includegraphics[scale=0.55]{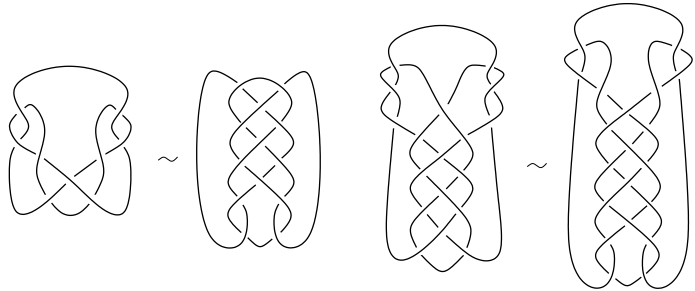}
  \caption{Two symmetric union presentations of the knot $8_8$
    with partial knots $4_1$ and $5_1$, respectively, and the knot 
    12a3 with distinct partial knots $C(3,1,2)$ and $C(2,5)$.}
  \label{fig:Knots-8_8-12a3}
  %% For each pair symmetric equivalence is obstructed by distinct partial knots.
  %% Two symmetric union presentations of the knot 11a36, 
  %% with distinct partial knots $C(3,1,2)$ and $C(2,5)$, and of 
  %% the knot 12a3, with distinct partial knots $C(3,4)$ and $C(2,6)$
\end{figure}

\begin{example}
  It is easy to see that symmetric Reidemeister moves do not change 
  the partial knots (see \textsection\ref{sub:PartialKnotInvariance}).  
  Consequently, if a knot $K$ can be represented by two symmetric unions 
  with distinct pairs of partial knots, then the two representations 
  cannot be equivalent under symmetric Reidemeister moves.  Two examples 
  of this type are depicted in Figure \ref{fig:Knots-8_8-12a3}. 
  The smallest known examples are the knots $8_8$ and $8_9$:
  for each of them we found two symmetric unions 
  with partial knots $4_1$ and $5_1$, respectively.
  This shows that partial knots need not be unique
  even for the class of two-bridge ribbon knots.
\end{example}

\begin{figure}[hbtp]
  \centering
  \includegraphics[scale=0.5]{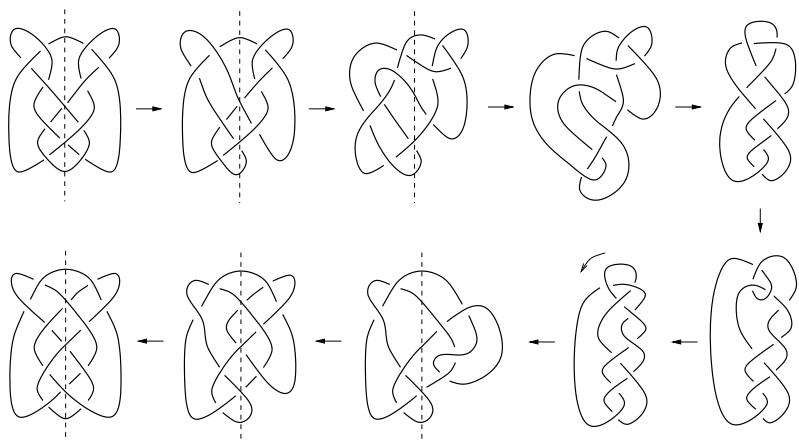}
  \caption{A transformation between two symmetric 
    union representations of the knot $9_{27}$.  
    The intermediate stages are not symmetric.}
  \label{fig:Transformation}
\end{figure}

Partial knots form an obvious obstruction for symmetric equivalence,
but there also exist examples of symmetric union representations
with the same partial knots, but which are not symmetrically equivalent.
Figure \ref{fig:Transformation} shows a transformation between 
the two symmetric union representations of the knot $9_{27}$ 
displayed in Figure \ref{fig:Knot-9-27} at the beginning of this article.

\begin{theorem}
  The two symmetric union diagrams shown in Figure \ref{fig:Knot-9-27} 
  both represent the knot $9_{27}$ and both have $5_2$ as partial knot.
  They are, however, not equivalent under symmetric Reidemeister moves
  as defined in \textsection\ref{sec:SymmetricReidemeister}.
  \qed
\end{theorem}

While some experimentation might convince you that 
this result is plausible, it is not so easy to prove.
We will give the proof in a forthcoming article
\cite{EisermannLamm:SymJones}, based on a two-variable 
refinement of the Jones polynomial for symmetric unions.
The diagrams displayed here are the first pair 
of an infinite family of two-bridge knots exhibited 
in \textsection\ref{sec:InfiniteFamily}.

\subsection{Knots with extra structure}

The study of symmetric diagrams and symmetric equivalence is meaningful 
also for other types of symmetries, or even more general constraints.
It can thus be seen as an instance of a very general principle,
which could be called \emph{knots with extra structure},  
and which seems worthwhile to be made explicit.

Generally speaking, we are given a class of diagrams satisfying some 
constraint and a set of (generalized) Reidemeister moves respecting
the constraint.  It is then a natural question to ask whether the 
equivalence classes under constrained moves are strictly smaller than 
those under usual Reidemeister moves (ignoring the constraint, 
e.g.\ breaking the symmetry).  If this is the case then two opposing
interpretations are possible:
\begin{itemize}
\item[(a)]
  We might have missed some natural but less obvious move that respects 
  the constraint.  Such a move should be included to complete our list.
\item[(b)]
  The constraint introduces some substantial obstructions 
  that cannot be easily circumvented.  
  The induced equivalence is an object in its own right.
\end{itemize}

In order to illustrate the point, let us cite some prominent examples,
which have developed out of certain quite natural constraints.

\begin{itemize}
\item
  Perhaps the most classical example of diagrams and moves
  under constraints is provided by alternating diagrams
  and Tait's flype moves, cf.\ \cite{MenascoThistlethwaite:1993}.
\item
  Braids form another important and intensely studied case.
  Here one considers link diagrams in the form of a closed braid
  and Markov moves, cf.\ \cite{Birman:1974,BirmanMenasco:2002}.
\end{itemize}

In these two settings the fundamental result is that constrained 
moves generate the same equivalence as unconstrained moves.  
In the following two examples, however, new classes 
of knots have emerged:

\begin{itemize}
\item
  Given a contact structure, one can consider knots that are 
  everywhere transverse (resp.\ tangent) to the plane field,
  thus defining the class of transverse (resp.\ legendrian) knots.  
  Again one can define equivalence by isotopies respecting 
  this constraint, and it is a natural question to what extent 
  this equivalence is a refinement of the usual equivalence,
  cf.\ \cite{BirmanMenasco:2006}.
\item
  Virtual knots can also be placed in this context: here one introduces 
  a new type of crossing, called virtual crossing, and allows suitably 
  generalized Reidemeister moves, cf.\ \cite{Kauffman:1999}.  
  Strictly speaking, this is an extension rather than a constraint,
  and classical knots inject into the larger class of virtual knots.
\end{itemize}

Considering symmetric unions, two nearby generalizations also seem promising:

\begin{itemize}
\item
  Analogous to diagrams that are symmetric with respect 
  to reflection, one can consider strongly amphichiral diagrams.
  Here the symmetry is a rotation of $180^\circ$ about a point, 
  which maps the diagram to itself reversing all crossings.
  Again there are some obvious moves respecting the symmetry,
  leading to a natural equivalence relation on the set
  of strongly amphichiral diagrams.
\item
  Since ribbon knots are in general not known 
  to be representable as symmetric unions, one could
  consider band presentations of ribbon knots and Reidemeister 
  moves respecting the band presentation.  The equivalence
  classes will thus correspond to ribbons modulo isotopy.
  For a given knot $K$ the existence and uniqueness questions 
  can be subsumed by asking how many ribbons are there for $K$.
\end{itemize}

Of course, the paradigm of ``knots with extra structure'' cannot be expected 
to produce any general answers; the questions are too diverse and often rather deep.  
Nevertheless, we think of it as a good generic starting point and a unifying 
perspective.  Its main merit is that it leads to interesting questions.
In the present article we will begin investigating 
this approach in the special case of symmetric unions.

%%%%%%%%%%%%%%%%%%%%%%%%%%%%%%%%%%%%%%%%%%%%%%%%%%%%%%%%%%%%%%%%%%%%%%%%%%%%%

\section{Symmetric diagrams and symmetric moves} 
\label{sec:SymmetricReidemeister}

Having seen some examples of symmetric unions that are equivalent 
by asymmetric Reidemeister moves, we wish to make precise what 
we mean by \emph{symmetric equivalence}.  As can be suspected,
this will be the equivalence relation generated by symmetric 
Reidemeister moves, but the details require some attention.

\subsection{Symmetric Reidemeister moves} \label{sub:SymmetricMoves}

We consider the euclidian plane $\R^2$ with the reflection
$\rho\colon \R^2 \to \R^2$, $(x,y) \mapsto (-x,y)$.
The map $\rho$ reverses orientation and its fix-point 
set is the vertical axis $\{0\} \times \R$.
A link diagram $D \subset \R^2$ is \emph{symmetric} 
with respect to this axis if and only if $\rho(D) = D$ 
except for crossings on the axis, which are necessarily reversed.

By convention we will not distinguish two symmetric diagrams 
$D$ and $D'$ if they differ only by an orientation preserving 
diffeomorphism $h \colon \R^2 \isoto \R^2$ respecting the symmetry, 
in the sense that $h \circ \rho = \rho \circ h$.

\begin{figure}[hbtp]
  \centering
  \includegraphics[scale=0.8]{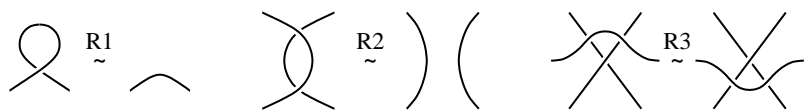}
  \caption{The classical Reidemeister moves (off the axis)}
  \label{fig:Rmoves}
\end{figure}

\begin{figure}[hbtp]
  \centering
  \includegraphics[scale=0.9]{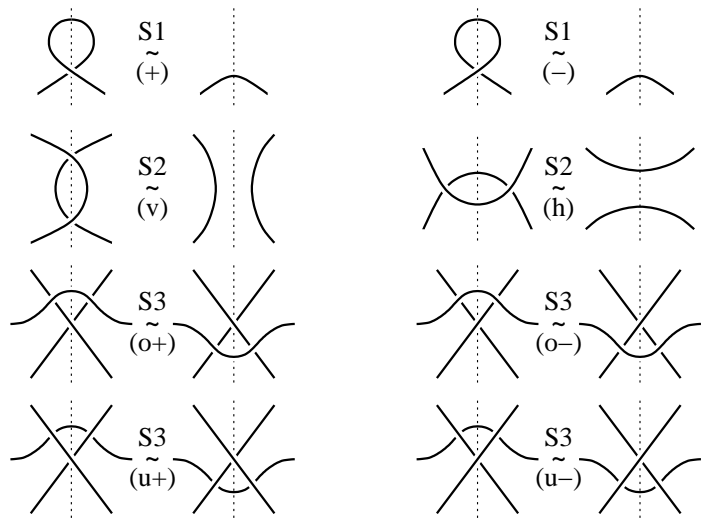}
  \caption{Symmetric Reidemeister moves on the axis}
  \label{fig:Smoves}
\end{figure}

\begin{figure}[hbtp]
  \centering
  \includegraphics[scale=0.8]{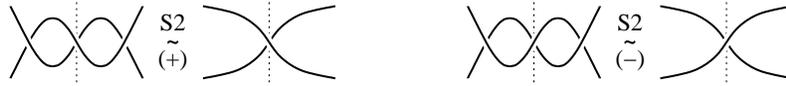}
  \caption{A symmetric move on two strands with three crossings}
  \label{fig:S2move}
\end{figure}

\begin{figure}[hbtp]
  \centering
  \includegraphics[scale=0.9]{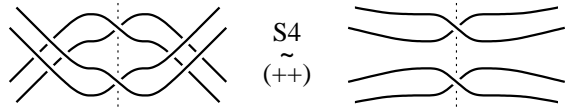}
  \caption{A symmetric move on four strands across the axis}
  \label{fig:S4move}
\end{figure}

\begin{definition}
  Given a symmetric diagram, a \emph{symmetric Reidemeister move} 
  with respect to the reflection $\rho$ is a move of the following type:
  \begin{itemize}
  \item
    A symmetric Reidemeister move off the axis, that is, an ordinary
    Reidemeister move, R1--R3 as depicted in Figure \ref{fig:Rmoves},
    carried out simultaneously with its mirror-symmetric 
    counterpart with respect to the reflection $\rho$.
  \item
    A symmetric Reidemeister move on the axis, 
    of type S1--S3 as depicted in Figure \ref{fig:Smoves},
    or a generalized Reidemeister move on the axis, 
    of type S2($\pm$) as depicted in Figure \ref{fig:S2move}, 
    or of type S4 as depicted in Figure \ref{fig:S4move}.
  \end{itemize}
\end{definition}

\vspace{-6pt}

\begin{remark}
  By convention the axis is not oriented in the local pictures, 
  so that we can turn Figures \ref{fig:Smoves}, 
  \ref{fig:S2move}, \ref{fig:S4move} upside-down.
  This adds one variant for each S1-, S2-, and S4-move shown here;
  the four S3-moves are invariant under this rotation.
  Moreover, the S4-move comes in four variants, obtained by 
  changing the over- and under-crossings on the axis.
\end{remark}

\subsection{Are these moves necessary?}

\begin{figure}[hbtp]
  \centering
  \includegraphics[scale=0.95]{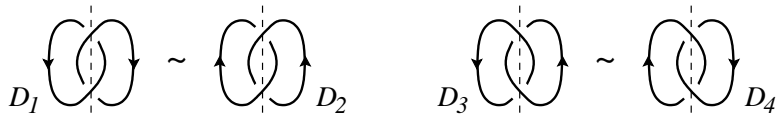}
  \caption{S2($\pm$)-moves and S4-moves are necessary.}
  \label{fig:S4needed}
\end{figure}

The emergence of the somewhat unusual moves S2($\pm$) 
and S4 may be surprising at first sight.  
One might wonder whether they are necessary 
or already generated by the other, simpler moves:

\begin{theorem}
  The four oriented link diagrams shown in Figure \ref{fig:S4needed}
  all represent the Hopf link with linking number $+1$.
  The pairs $D_1 \sim D_2$ and $D_3 \sim D_4$ are equivalent 
  via symmetric Reidemeister moves, but $D_1,D_2$ 
  are only asymmetrically equivalent to $D_3,D_4$.
  Moreover, the symmetric equivalence $D_1 \sim D_2$ 
  cannot be established without using S2($\pm$)-moves, 
  and the symmetric equivalence $D_3 \sim D_4$ 
  cannot be established without using S4-moves.
\end{theorem}

\begin{proof}
  \newcommand{\pic}[1]{\raisebox{-0.9ex}{\includegraphics[height=3ex]{#1}}}
  The symmetric equivalences $D_1 \sim D_2$ 
  and $D_3 \sim D_4$  are easily established,
  and will be left as an amusing exercise.
  Less obvious is the necessity of moves S2($\pm$) and S4.
  Likewise, $D_2$ and $D_3$ are asymmetrically equivalent,
  but we need an obstruction to show that they cannot be symmetrically equivalent.
  Given an oriented diagram $D$, we consider the points 
  on the axis where two distinct components cross.
  To each such crossing we associate an element in 
  the free group $F = \langle s,t,u,v \rangle$ as follows:
  \begin{xalignat*}{4}
    \pic{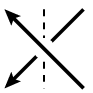}  & \mapsto s^{+1} &
    \pic{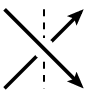} & \mapsto t^{+1} &
    \pic{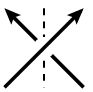}    & \mapsto u^{+1} &
    \pic{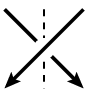}  & \mapsto v^{+1} \\
    \pic{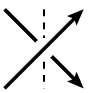} & \mapsto s^{-1} &
    \pic{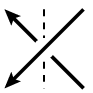}  & \mapsto t^{-1} &  
    \pic{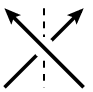}    & \mapsto u^{-1} &
    \pic{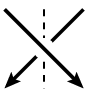}  & \mapsto v^{-1} 
  \end{xalignat*}
  
  Traversing the axis from top to bottom we read a word on the alphabet
  $\{ s^\pm, t^\pm, u^\pm, v^\pm \}$, which defines an element $w(D) \in F$.
  It is an easy matter to verify how symmetric Reidemeister moves 
  affect $w(D)$.  Moves off the axis have no influence.
  S1-moves are neglected by construction.
  S2(v)-moves change the word by a trivial relation
  so that $w(D) \in F$ remains unchanged.
  S2(h)-moves and S3-moves have no influence.  
  An S2($\pm$)-move can change one factor $u \leftrightarrow v$,
  but leaves factors $s$ and $t$ unchanged.
  An S4-move, finally, interchanges two adjacent factors.

  In our example we have $w(D_1) = u^2$ and $w(D_2) = v^2$,
  so at least two S2($\pm$)-moves are necessary in the transformation.
  Furthermore, $w(D_3) = s t$ and $w(D_4) = t s$, so at least
  one S4-move is necessary in the transformation.
  Finally, no symmetric transformation can change 
  $D_1$ or $D_2$ into $D_3$ or $D_4$.
\end{proof}

\begin{remark}[orientation]
  One might object that the preceding proof introduces 
  the orientation of strands as an artificial subtlety.
  Denoting for each oriented diagram $D$ the underlying 
  unoriented diagram by $\bar{D}$, we see that 
  $\bar{D}_1 = \bar{D}_2$ and $\bar{D}_3 = \bar{D}_4$ 
  are identical unoriented diagrams. 
  Orientations obviously simplify the argument, 
  but it is worth noting that the phenomenon 
  persists for unoriented knot diagrams as well:
\end{remark}

\begin{corollary}
  The unoriented diagrams $\bar{D}_2$ and $\bar{D}_3$ 
  are not symmetrically equivalent.
\end{corollary}

\begin{proof}
  If $\bar{D}_2$ were symmetrically equivalent to $\bar{D}_3$,
  then we could equip $\bar{D}_2$ with an orientation, say $D_2$,
  and carry it along the transformation to end up with some
  orientation for $\bar{D}_3$.  Since the linking number
  must be $+1$, we necessarily obtain $D_3$ or $D_4$.
  But $w(D_2) = v^2$ can be transformed neither into $w(D_3) = s t$ 
  nor $w(D_4) = t s$. This is a contradiction.
\end{proof}

\begin{corollary}
  The moves S2$(\pm)$ and S4 are also necessary for
  the symmetric equivalence of unoriented diagrams.
\end{corollary}

\begin{proof}
  The following trick allows us to apply
  the above argument to unoriented links:
  we take the diagrams of Figure \ref{fig:S4needed}
  and tie a non-invertible knot into each component, 
  symmetrically on the left and on the right.
  This introduces an intrinsic orientation.
  %% The details are more complicated, 
  %% but the argument is essentially the same.
\end{proof}

\begin{remark}[linking numbers]
  The proof of the theorem shows that the composition 
  $\bar{w} \colon \{\begin{smallmatrix} \text{oriented} \\ 
    \text{diagrams}\end{smallmatrix}\} \to F \to \Z^3$ defined by
  $s \mapsto (1,0,0)$, $t \mapsto (0,1,0)$, $u,v \mapsto (0,0,1)$
  is invariant under \emph{all} symmetric Reidemeister moves.
  For example $\bar{w}(D_1) = \bar{w}(D_2) = (0,0,2)$ and 
  $\bar{w}(D_3) = \bar{w}(D_4) = (1,1,0)$ yields the obstruction
  to symmetric equivalence used above.  The invariant $\bar{w}$ can be 
  interpreted as a refined \emph{linking number} for crossings on the axis.  
  This already indicates that the symmetry constraint
  may have surprising consequences.
\end{remark}

\begin{remark}[symmetric unions]
  While refined linking numbers may be useful for symmetric diagrams in general, 
  such invariants become useless when applied to \emph{symmetric unions}, 
  which are our main interest.  In this more restrictive setting 
  we only have \emph{one} component.  When trying to imitate
  the above construction, S1-moves force the relation $s = t = 1$.  
  Moreover, orientations are such that a crossing on the axis 
  always points ``left'' or ``right'' but never ``up'' or ``down'',
  so factors $u^\pm$ and $v^\pm$ never occur.
\end{remark}

\subsection{Invariance of partial knots} \label{sub:PartialKnotInvariance}

Recall that for every symmetric union diagram $D$ 
we can define partial diagrams $D_-$ and $D_+$ as follows: 
first, we resolve each crossing on the axis by cutting it open 
according to $\aocr \mapsto \avcr$ or $\aucr \mapsto \avcr$.
The result is a diagram $\hat{D}$ without any
crossings on the axis.  If we suppose that $D$ 
is a symmetric union, then $\hat{D}$ is a connected sum, 
which can then be split by a final cut $\ahcr \mapsto \avcr$.  
We thus obtain two disjoint diagrams: 
$D_-$ in the half-space $H_- = \{ (x,y) \mid x<0 \}$,
and $D_+$ in the half-space $H_+ = \{ (x,y) \mid x>0 \}$.
The knots $K_-$ and $K_+$ represented by $D_-$ and $D_+$,
respectively, are called the \emph{partial knots} of $D$.
Since $D$ was assumed symmetric, $K_+$ and $K_-$ 
are mirror images of each other.

\begin{proposition}
  For every symmetric union diagram $D$ the partial knots $K_-$ and $K_+$ 
  are invariant under symmetric Reidemeister moves.
\end{proposition}

\begin{proof}
  This is easily seen by a straightforward case-by-case verification.
\end{proof}

\subsection{Horizontal and vertical flypes}

The symmetric Reidemeister moves displayed above give 
a satisfactory answer to the local equivalence question.
There are also some semi-local moves that merit attention,
most notably flype moves.

\begin{figure}[hbtp]
  \centering
  \includegraphics[scale=0.8]{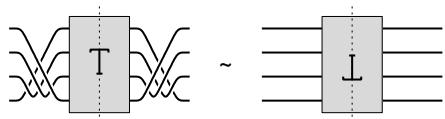}
  \caption{A horizontal flype (across the axis)}
  \label{fig:FlypeAcrossAxis}
\end{figure}

\begin{proposition}
  Every horizontal flype across the axis, 
  as depicted in Figure \ref{fig:FlypeAcrossAxis},
  can be decomposed into a finite sequence 
  of symmetric Reidemeister moves.
  \qed
\end{proposition}

\begin{figure}[hbtp]
  \centering
  \includegraphics[scale=0.8]{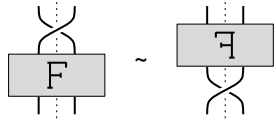}
  \caption{A vertical flype (along the axis)}
  \label{fig:FlypeAlongAxis}
\end{figure}

\begin{definition}
  A \emph{vertical flype along the axis} is a move 
  as depicted in Figure \ref{fig:FlypeAlongAxis}, where
  the tangle $F$ can contain an arbitrary diagram 
  that is symmetric with respect to the axis.
\end{definition}

\begin{figure}[hbtp]
  \centering
  \includegraphics[scale=0.8]{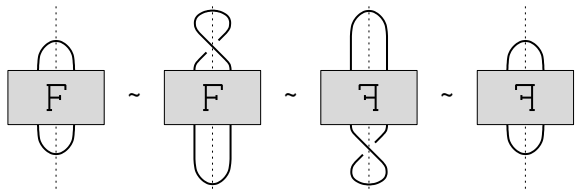}
  \caption{A flype rotating the entire diagram}
  \label{fig:FlypeAll}
\end{figure}

\begin{example}
  Strictly speaking a flype is not a local move, because 
  the tangle $F$ can contain an arbitrarily complicated diagram.
  Such a flype allows us, for example, to realize a rotation 
  of the entire diagram around the axis, as depicted
  in Figure \ref{fig:FlypeAll}.
\end{example}

While a horizontal flype can be achieved by symmetric 
Reidemeister moves, this is in general not possible for
a vertical flype: when decomposed into Reidemeister moves,
the intermediate stages are in general no longer symmetric.
This is also manifested in the following observation:

\begin{proposition}
  A vertical flype changes the partial knots in a well-controlled way, 
  from $K_- \cs L_-$ and $K_+ \cs L_+$ to $K_- \cs L_+$ and $K_+ \cs L_-$,
  where $(K_-,K_+)$ and $(L_-,L_+)$ are pairs of mirror images.
  In general this cannot be realized by symmetric Reidemeister moves.
  \qed
\end{proposition}

\subsection{Connected sum}

As a test-case for symmetric equivalence,
we wish to construct a connected sum for symmetric unions
and show that it shares some properties with the usual connected sum.
This is by no means obvious, and the first problem will be the 
very definition: is the connected sum well-defined on equivalence classes?  
The fact that the answer is affirmative can be seen as 
a confirmation of our chosen set of Reidemeister moves.

In order to define a connected sum of diagrams
we have to specify which strands will be joined.
To this end we consider pointed diagrams as follows.  

\begin{definition}
  Each symmetric union diagram $D$ traverses the axis 
  at exactly two points that are not crossings.  
  We mark one of them as the basepoint of $D$.
  The result will be called a \emph{pointed diagram}.
  Given two symmetric union diagrams $D$ and $D'$
  that are pointed and oriented, we can define their 
  \emph{connected sum} $D \cs D'$ as indicated
  in Figure \ref{fig:ConnectedSum}.  
  The result is again a symmetric union diagram 
  that is pointed and oriented.

  \begin{figure}[hbtp]
    \centering
    \includegraphics[scale=0.8]{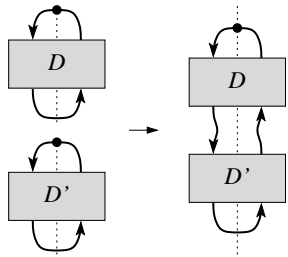}
    \caption{Connected sum $D \cs D'$ of two symmetric union diagrams 
      $D$ and $D'$, each equipped with a basepoint and an orientation}
    \label{fig:ConnectedSum}
  \end{figure} 
  
  More explicitly, we start with the distant union $D \sqcup D'$ 
  by putting $D$ above $D'$ along the axis.  Both diagrams intersect
  the axis transversally at two points each.  We choose the
  unmarked traversal of $D$ and the marked traversal of $D'$ 
  and join the strands to form the connected sum.  
  If other strands (with crossings on the axis) are 
  between them, we pass over all of them by convention.  
  If the orientations do not match, we perform an S1$+$ 
  move on one of the strands before joining them.
\end{definition}

All symmetric moves discussed previously generalize 
to pointed diagrams: in each instance the basepoint 
is transported in the obvious way. %% [FIGURE?]
The upshot is the following result:

\begin{theorem}
  The connected sum induces a well-defined operation on equivalence 
  classes modulo symmetric Reidemeister moves and flypes.
  More explicitly, this means that $D_1 \sim D_2$ and 
  $D'_1 \sim D'_2$ imply $D_1 \cs D'_1 \sim D_2 \cs D'_2$.
  The connected sum operation is associative and has the class 
  of the trivial diagram as two-sided unit element.
\end{theorem}

\begin{proof}
  Consider a symmetric Reidemeister move performed 
  on the diagram $D'$.  If the basepoint is not concerned 
  then the same move can be carried out in $D \cs D'$.
  Only two special cases need clarification: an S1-move 
  (or more generally a flype move) on $D'$ affecting 
  the basepoint translates into a flype move on $D \cs D'$.
  An S3-move on $D'$ affecting the basepoint can be translated 
  to a sequence of symmetric Reidemeister moves on $D \cs D'$.
  The situation is analogous for $D$ concerning the 
  unmarked traversal of the axis; the verifications 
  are straightforward.
\end{proof}

\subsection{Open questions}

The connected sum of symmetric unions, as defined above,
is associative but presumably not commutative.
The usual trick is to shrink $D'$ and to slide it along $D$ 
so as to move from $D \cs D'$ to $D' \cs D$, but this 
transformation is not compatible with our symmetry constraint.  
Even though non-commutativity is a plausible consequence, 
this does not seem easy to prove.  

\begin{question}
  Is the connected sum operation on symmetric unions
  non-commutative, as it seems plausible? 
  How can we prove it?  Does this mean that
  we have missed some less obvious but natural move?
  Or is it an essential feature of symmetric unions?
\end{question}

\begin{remark}
  On the one hand non-commutativity may come as a surprise
  for a connected sum operation of knots.  On the other hand, 
  the connected sum of symmetric unions is halfway between 
  knots and two-string tangles, and the latter are highly
  non-commutative.  The theory of symmetric unions retains
  some of this two-string behaviour.

  Although only loosely related, we should also like to point out
  that similar phenomena appear for virtual knots.
  There the connected sum is well-defined only for long knots,
  corresponding to a suitable marking how to join strands.
  Moreover,  the connected sum for long virtual knots
  is not commutative \cite{SilverWilliams:2006}.  
\end{remark}

\begin{question}
  %% What is the centre, i.e.\ which elements commute with all others?
  Symmetric unions of the form $K_+ \cs K_-$ belong to the centre:
  the usual trick of shrinking and sliding $K_\pm$ along the strand
  still works in the symmetric setting.  Are there any other central elements?
\end{question}

\begin{question}
  \newcommand{\pic}[1]{\raisebox{-0.9ex}{\includegraphics[height=3ex]{#1}}}
  What are the invertible elements, satisfying $D \cs D' = \pic{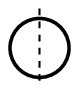}$?
  An invertible symmetric union diagram necessarily represents the unknot.  
  It is not clear, however, if it is equivalent to the unknot by symmetric moves.
\end{question}

\begin{question}
  Do we have unique decomposition into prime elements?
\end{question}

\begin{question}
  Some geometric invariants such as bridge index,
  braid index, and genus, can be generalized to the setting
  of symmetric unions, leading to the corresponding
  notions of symmetric bridge index, symmetric braid index,
  and symmetric genus.  Do they have similar properties
  as in the classical case, i.e.\  is the unknot detected
  and does connected sum translate into addition?
\end{question}

\begin{remark}
  If we had some additive invariant $\nu$ and a non-trivial
  symmetric union representation $U$ of the unknot with $\nu(U) > 0$,
  then every symmetric union diagram $D$ would yield an infinite family
  $D \cs U^{\cs k}$ of distinct diagrams representing the same knot.
\end{remark}

%%%%%%%%%%%%%%%%%%%%%%%%%%%%%%%%%%%%%%%%%%%%%%%%%%%%%%%%%%%%%%%%%%%%%%%%%%%%%

\section{Inequivalent symmetric union representations} \label{sec:InfiniteFamily}

\subsection{An infinite family}

In this section we exhibit an infinite family of symmetric unions 
which extend the phenomenon observed for the diagrams of $9_{27}$.
Notice that we will be dealing with prime knots, so this non-uniqueness 
phenomenon is essentially different from the non-uniqueness caused 
by the non-commutativity of the connected sum operation.

\begin{figure}[hbtp]
  \centering
  \includegraphics[scale=0.55]{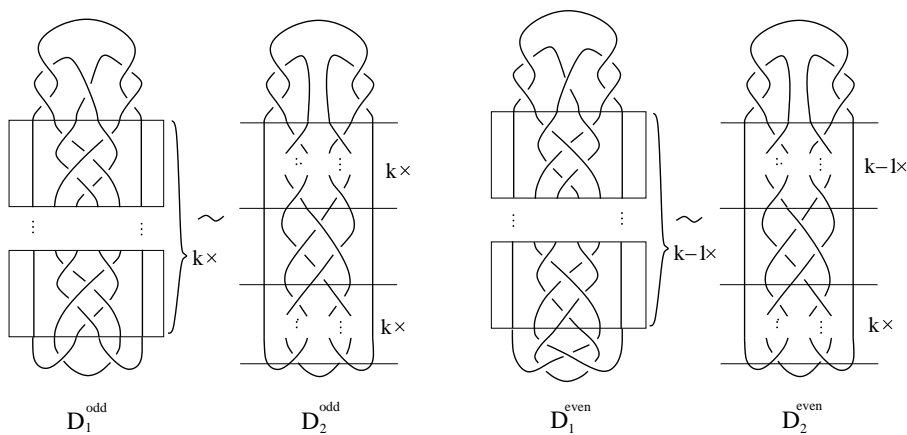}
  \caption{Statement of the theorem}
  \label{family:statement}
\end{figure}

\begin{definition}
  For each integer $n \ge 2$ we define two symmetric union diagrams 
  $D_1(n)$ and $D_2(n)$ as follows.  We begin with the connected sum
  $C(2,n) \sharp C(2,n)^*$ and insert crossings on the axis 
  as indicated in Fig.\ \ref{family:statement},
  distinguishing the odd case $n=2k+1$ and the even case $n=2k$.  
\end{definition}

\begin{theorem}
  For each $n \ge 2$ the diagrams $D_1(n)$ and $D_2(n)$ 
  can be transformed one into another by a sequence 
  of Reidemeister moves, not respecting the symmetry:
  \[
  D_1(n) \sim D_2(n) \sim \begin{cases}
    S\left( (2n+1)^2, 2n^2 \right) & \text{if $n$ is odd} , \\
    S\left( (2n+1)^2, 2n^2-1 \right) & \text{if $n$ is even} .
  \end{cases}
  \]
\end{theorem}

Here $S(p,q)$ is Schubert's notation for two-bridge knots, 
see \cite[\textsection2.1]{Kawauchi:1996}.  

\begin{example}
  For $n=2$ we obtain two mirror-symmetric diagrams 
  $D_1(2)$ and $D_2(2)$ of the knot $8_9$,
  which turn out to be symmetrically equivalent.
  For $n=3$ we obtain the two symmetric union representations 
  of $9_{27}$ depicted in Fig.\ \ref{fig:Knot-9-27}.

  These and the following cases yield two symmetric union representations 
  of the two-bridge knots $K(a,b)=C(2a,2,2b,-2,-2a,2b)$ with $b=\pm1$, 
  up to mirror images: more explicitly, we so obtain the knots 
  $8_9    = K(-1,-1)$ for $n=2$, 
  $9_{27} = K(-1,1)$  for $n=3$, 
  $10_{42}= K(1,1)$   for $n=4$, 
  $11a96  = K(1,-1)$  for $n=5$,
  $12a715 = K(-2,-1)$ for $n=6$, 
  $13a2836= K(-2,1)$  for $n=7$. 
  They all have genus $3$ and their crossing number is $6+n$.
\end{example}

After some experimentation you might find it plausible 
that $D_1(n)$ and $D_2(n)$ are not \emph{symmetrically} 
equivalent for $n \ge 3$.  
Notice, however, that the obvious obstruction fails:
by construction, both have the same partial knots $C(2,n)$ and $C(2,n)^*$.
Their non-equivalence will be studied in \cite{EisermannLamm:SymJones}
where we develop the necessary tools.

\begin{proof} % [Proof of the theorem]
  We first analyze the braid $\beta_1$ that is 
  shown boxed in diagram $D^\text{odd}_{1}$. 
  Using the braid relations we have 
  \begin{align*}
    \beta_1
    & = \sigma_2^{-1}\sigma_4\sigma_3^{-1}\sigma_2^{-1}\sigma_4\sigma_3 
    = \sigma_4\sigma_2^{-1}\sigma_3^{-1}\sigma_2^{-1}\sigma_4\sigma_3 
    \\
    & = \sigma_4\sigma_3^{-1}\sigma_2^{-1}\sigma_3^{-1}\sigma_4\sigma_3 
    = \sigma_4\sigma_3^{-1}\sigma_2^{-1}\sigma_4\sigma_3\sigma_4^{-1} 
  \end{align*}
  Therefore $\beta_1^k = \sigma_4 \sigma_3^{-1}\sigma_2^{-k}\sigma_4^{k}\sigma_3\sigma_4^{-1}$.
  For the braid $\beta_2$, shown boxed in diagram $D^\text{even}_{1}$, we have similarly
  \begin{equation*}
  \beta_2 = \sigma_2^{-1}\sigma_4\sigma_3\sigma_2^{-1}\sigma_4\sigma_3^{-1}
          = \sigma_2^{-1}\sigma_3\sigma_2^{-1}\sigma_4\sigma_3^{-1}\sigma_2
  \end{equation*}
  and $\beta_2^k = \sigma_2^{-1}\sigma_3\sigma_2^{-k}\sigma_4^k\sigma_3^{-1}\sigma_2$.
  With this information at hand, we pursue the odd and even cases separately.
  
  \medskip\noindent\textit{First case: $n$ is odd.}
  The simplification of $D^\text{odd}_{1}$ done by computing $\beta_1^k$ 
  is shown in diagram $D^\text{odd}_{1'}$ in Fig.\ \ref{family:proof:odd}.
  This diagram can be further transformed, yielding diagram $D^\text{odd}_{1''}$ 
  which   is in two-bridge form.  Its Conway notation is $C(2,k,2,1,2,-k-1)$. 
  
  Diagram $D^\text{odd}_2$ in Fig.\ \ref{family:statement} 
  simplifies to $D^\text{odd}_{2'}$ in Fig.\ \ref{family:proof:odd} 
  because certain crossings are cancelled.  Further transformation 
  gives its two-bridge form, shown in diagram $D^\text{odd}_{2''}$.
  Its Conway notation is $C(2,k,2,2,-2,-k)$.  
  The continued fractions for both knots evaluate  
  to $\frac{(4k+3)^2}{8k^2+8k+2} = \frac{(2n+1)^2}{2n^2}$, 
  so both knots are equal.  

  \begin{figure}[hbtp]
    \centering
    \includegraphics[scale=0.6]{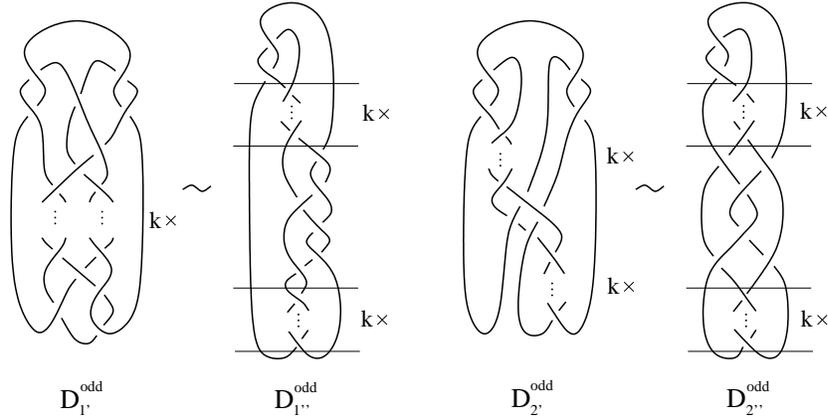}
    \caption{Proof in the odd case}
    \label{family:proof:odd}
  \end{figure}

  \medskip\noindent\textit{Second case: $n$ is even.}
  We simplify the braid $\beta_2^{k-1}\sigma_2^{-1}\sigma_4\sigma_3\sigma_2^{-1}\sigma_4$
  occurring in diagram $D^\text{even}_{1}$ of Fig.\ \ref{family:statement}:
  using the formula for $\beta_2^{k-1}$ and applying braid relations we get 
  \[
  \sigma_2^{-1}\sigma_3(\sigma_2^{-1}\sigma_4)^{k-1}\sigma_4\sigma_3\sigma_2^{-1}
  \]
  which is depicted in diagram $D^\text{even}_{1'}$ 
  in Fig.\ \ref{family:proof:even}. 

  \begin{figure}[hbtp]
    \centering
    \includegraphics[scale=0.6]{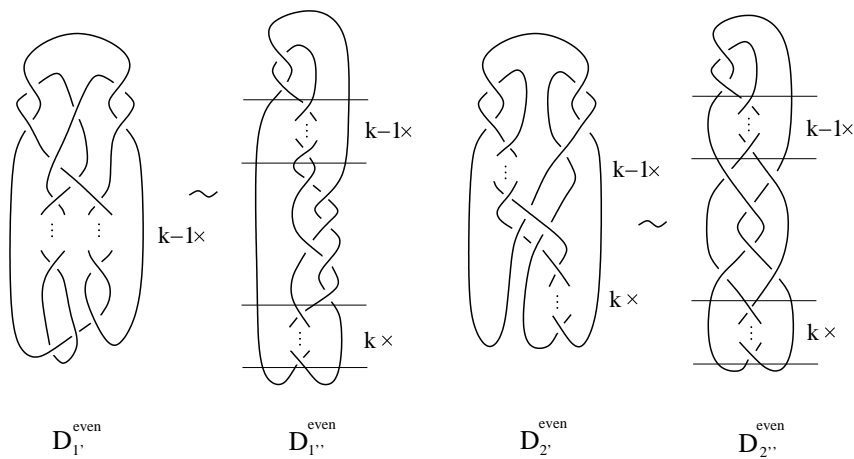}
    \caption{Proof in the even case}
    \label{family:proof:even}
  \end{figure}

  The transformation to two-bridge form is similar to the odd case and 
  we get the knot $C(2,k,-2,-1,-2,-k)$ shown in diagram $D^\text{even}_{1''}$. 
  The simplification of diagram $D^\text{even}_{2}$ in Fig.\ \ref{family:statement} 
  to diagram $D^\text{even}_{2'}$ in Fig.\ \ref{family:proof:even} 
  is straightforward, and diagram $D^\text{even}_{2''}$ allows us 
  to read off its two-bridge form $C(2,k-1,2,-2,-2,-k)$.  
  The continued fractions for both knots evaluate to 
  $\frac{(4k+1)^2}{8k^2-1} = \frac{(2n+1)^2}{2n^2-1}$. 
\end{proof}

\subsection{Open questions}

As we have seen, certain ribbon knots have 
more than one symmetric representation.
We have not succeeded in finding such
an ambiguity for the two smallest ribbon knots:

\begin{question}
  Can the unknot be represented by symmetric union diagrams 
  belonging to more than one equivalence class?  It is known
  that the partial knots of a symmetric union representation of 
  the unknot are necessarily trivial, see \cite[Theorem 3.5]{Lamm:2000}.
\end{question}

\begin{question}
  Can the knot $6_1$ be represented by symmetric union diagrams 
  belonging to more than one equivalence class?
\end{question}

\begin{question}
  Is the number of equivalence classes of symmetric unions
  representing a given knot $K$ always finite?
  Does non-uniqueness have some geometric meaning?
  For example, do the associated ribbon bands 
  differ in some essential way?
\end{question}

%%%%%%%%%%%%%%%%%%%%%%%%%%%%%%%%%%%%%%%%%%%%%%%%%%%%%%%%%%%%%%%%%%%%%%%%%%%%%

\bibliographystyle{plain}
\bibliography{symequiv}

\begin{thebibliography}{10}

\bibitem{Birman:1974}
J.~S. Birman.
\newblock {\em Braids, links, and mapping class groups}.
\newblock Princeton University Press, Princeton, N.J., 1974.

\bibitem{BirmanMenasco:2002}
J.~S. Birman and W.~W. Menasco.
\newblock On {M}arkov's theorem.
\newblock {\em J. Knot Theory Ramifications}, 11(3):295--310, 2002.

\bibitem{BirmanMenasco:2006}
J.~S. Birman and W.~W. Menasco.
\newblock Stabilization in the braid groups. {II}. {T}ransversal simplicity of
  knots.
\newblock {\em Geom. Topol.}, 10:1425--1452, 2006.

\bibitem{CassonGordon:1975}
A.~J. Casson and C.~McA. Gordon.
\newblock Cobordism of classical knots.
\newblock In {\em \`A la recherche de la topologie perdue}, volume~62 of {\em
  Progr. Math.}, pages 181--199. Birkh\"auser, Boston, MA, 1986.

\bibitem{ChaLivingston:2003v6}
J.~C. Cha and C.~Livingston.
\newblock Unknown values in the table of knots.
\newblock arxiv:math.GT/0503125v6, 2006.

\bibitem{EisermannLamm:SymJones}
M.~Eisermann and C.~Lamm.
\newblock A two-variable refinement of the {J}ones polynomial for symmetric
  unions of knots, 2007.
\newblock In preparation.

\bibitem{Knotscape}
J.~Hoste and M.~Thistlethwaite.
\newblock {K}not{S}cape -- providing convenient access to tables of knots.
\newblock http://www.math.utk.edu/{\textasciitilde}morwen/knotscape.html.

\bibitem{Kauffman:OnKnots}
L.~H. Kauffman.
\newblock {\em On knots}, volume 115 of {\em Annals of Mathematics Studies}.
\newblock Princeton University Press, Princeton, NJ, 1987.

\bibitem{Kauffman:1999}
L.~H. Kauffman.
\newblock Virtual knot theory.
\newblock {\em European J. Combin.}, 20(7):663--690, 1999.

\bibitem{Kawauchi:1996}
A.~Kawauchi.
\newblock {\em A survey of knot theory}.
\newblock Birkh\"auser Verlag, Basel, 1996.

\bibitem{KinoshitaTerasaka:1957}
S.~Kinoshita and H.~Terasaka.
\newblock On unions of knots.
\newblock {\em Osaka Math. J.}, 9:131--153, 1957.

\bibitem{Lamm:2000}
C.~Lamm.
\newblock Symmetric unions and ribbon knots.
\newblock {\em Osaka J. Math.}, 37(3):537--550, 2000.

\bibitem{Lamm:Talk}
C.~Lamm.
\newblock Symmetric union presentations for 2-bridge ribbon knots.
\newblock arxiv:math.GT/0602395, 2006.

\bibitem{Lisca:2007}
P.~Lisca.
\newblock Lens spaces, rational balls and the ribbon conjecture.
\newblock {\em Geom. Topol.}, 11:429--472, 2007.

\bibitem{Livingston:2005}
C.~Livingston.
\newblock A survey of classical knot concordance.
\newblock In {\em Handbook of knot theory}, pages 319--347. Elsevier B. V.,
  Amsterdam, 2005.

\bibitem{MenascoThistlethwaite:1992}
W.~W. Menasco and M.~Thistlethwaite.
\newblock The classification of alternating links.
\newblock {\em Ann. of Math. (2)}, 138(1):113--171, 1993.

\bibitem{Perko:1974}
K.~A. Perko, Jr.
\newblock On the classification of knots.
\newblock {\em Proc. Amer. Math. Soc.}, 45:262--266, 1974.

\bibitem{Rolfsen:1990}
D.~Rolfsen.
\newblock {\em Knots and links}, volume~7 of {\em Mathematics Lecture Series}.
\newblock Publish or Perish Inc., Houston, TX, 1990.
\newblock Corrected reprint of the 1976 original.

\bibitem{SilverWilliams:2006}
D.~S. Silver and S.~G. Williams.
\newblock Alexander groups of long virtual knots.
\newblock {\em J. Knot Theory Ramifications}, 15(1):43--52, 2006.

\end{thebibliography}

%%%%%%%%%%%%%%%%%%%%%%%%%%%%%%%%%%%%%%%%%%%%%%%%%%%%%%%%%%%%%%%%%%%%%%%%%%%%%
\end{document}